\documentclass[12pt]{article}
\usepackage{graphicx}
\usepackage{amssymb}
\usepackage{enumerate}
\usepackage{url}
\newcounter{fig}

\title{Marching in squares}
\author{Burkard Polster and Marty Ross}
\begin{document}

\maketitle

In the following we'll consider very mathematical (if not very plausible) geometric marching. Our marchers will exhibit  beautiful mathematics, some familiar and some less so.  In a final summary section we discuss the point of all the fancy footwork.

\section{Two squares make a square, ... or not} It is pleasing to know that there are parts of the world where everyone loves squares. What a paradise it must be, a country where people are forever marching in perfect square formation.\footnote{Google ``Japanese precision marching'' and ``Chinese army marching" to be amazed by what is possible in this respect.}

\begin{figure}[h]
\centerline{\includegraphics[scale=.7]{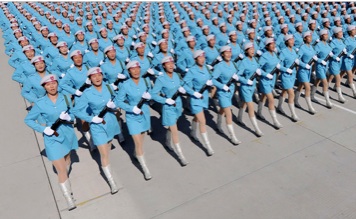}}
\end{figure}

Well, maybe not. Geometric marching is probably not as much fun when its main purpose is to please the Glorious Leader.
Still, those marching squares are impressive. And, given the unfortunate folk will be marching anyway, we have a great idea for a very mathematical flourish.

Our plan is to have two identical squares of marchers, each square performing the usual stunning steps. Then, the grand finale will consist of the squares being merged into one big (Red) square.

\begin{figure}[h]
\centerline{\includegraphics[scale=.5]{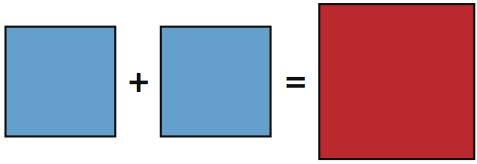}}
\end{figure}

It'd definitely be a showstopper but first there are details to sort out. The interlacing of the squares will be tricky, requiring planning, practice and fancy footwork.

We also have to decide the size of squares to use, which would need to consist of a suitable number of marchers. For example, beginning with two tiny $2 \times 2$ squares wouldn't work: that would give us 8 marchers in total, one short of the number required to rearrange into a $3 \times 3$ square. Similarly, beginning with two $3 \times 3$ squares of marchers would mean we have 18 marchers, too many to form a $4 \times 4$ square and insufficient for a $5 \times 5$ square.

Hmmm. This will take some figuring, but we should be able to do it. We'll try $4 \times 4$, then $5 \times 5$ and so on, and eventually we should have in hand the smallest squares that work.

For now, let's leave that calculation and just assume we've located the smallest squares that can be merged. Then, imagining we have sufficient marchers to occupy those squares, we can plan the marching steps.

Let's begin with an empty red quadrangle of just the right size to accommodate all our marchers. Then, a stylish approach would be to have the two identical squares of marchers enter the quadrangle from opposite sides.

\begin{figure}[h]
\centerline{\includegraphics[scale=.5]{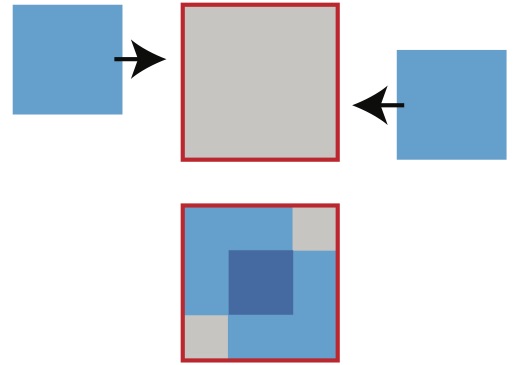}}
\end{figure}

At this stage the little grey squares are unoccupied and the blue squares are overlapping. It would be crowded in the middle dark blue square, but that's not a problem: we can simply arrange for half of the marchers to stand on the shoulders of the others.

Finally, we'll have the marchers leap off their comrades' shoulders and into the empty grey squares, a spectacular finish to our merging of the two blue squares into the Red square. Ta da!

But wait a minute. If, as planned, we have precisely the correct number of leaping marchers to fill the grey squares then this can be represented by the following picture.

\begin{figure}[h]
\centerline{\includegraphics[scale=.5]{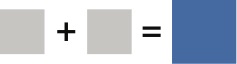}}
\end{figure}

Uh oh. We assumed that we began with the very smallest squares that could be merged to make a larger square. Yet, somehow we created even smaller squares that would work. How can that be?

Simply, it cannot be. For our smallest possible squares to result in even smaller squares that work is a plain logical impossibility. The unavoidable conclusion is that no squares can be merged in the way we had originally contemplated. So much for our plans to impress the Glorious Leader with a great marching finale.  However, perhaps he'll be impressed by some intriguing mathematics that emerges from our failed attempt. 

What  we have outlined above is a {\em proof by contradiction}. We began by assuming that certain squares were possible, and that assumption resulted in a   logical impossibility. This contradiction proves that our original assumption was wrong, and that no such squares can exist. 

We can now reconsider this conclusion in terms of numbers rather than squares. What we have demonstrated is that there are no natural numbers  $m$ and $n$  that solve the equation
$$2m^2 = n^2\,.$$
Rearranging, it follows that there is no fraction $n/m$ that solves the equation
 $$(n/m)^2 = 2\,.$$
So, there is no fraction whose square is 2, which is exactly the same as saying $\sqrt 2$ is not a fraction. That is, as you may have realised a while back, our marching ponderings prove that $\sqrt 2$ is irrational.

There are many proofs of the irrationality of $\sqrt 2$, however  the one underlying our marching square story above is probably our favourite. We learned it from the great John Conway, who attributes it to mathematician Stanley Tennenbaum; see  \cite{conway} and \cite{conway1}.

\section{Time to cheat}

Now all this is well and good but  we suspect that our mathematical ponderings won't impress our square-loving Glorious Leader. So, instead of completely giving up on the idea of merging squares into larger squares, we'll bend  the rules a little. 

Another way of expressing our square merging impossibility is to say that there is no {\em Pythagorean triple} of the form $m^2+m^2=n^2$. However, the Pythagorean triple $3^2+4^2=5^2$ demonstrates that we can get close: there is at least one  Pythagorean triple of the form $m^2+(m+1)^2=n^2$. Now if there were such a triple of numbers  for  reasonably large $m$ and $n$, then this would mean that basically indistinguishable squares containing $m^2$ and $(m+1)^2$ marchers could be merged into a larger square. And, hopefully the Glorious Leader wouldn't notice our little cheat.

 Just trying some numbers we find the solution $20^2+21^2=29^2$. To find more solutions we look to solve the equation 
$$m^2 + (m+1)^2=n^2\,.$$
Clearly the left side of this equation is always an odd number and so $n$ must be odd as well. Setting $n=2r+1$, we then have the 
equation
$$2m^2+2m+1=4r^2+4r+1\, ,$$ 
or, equivalently, 
$$r(r+1)+r(r+1)=m(m+1)\,.$$

That's very interesting and suggests a second way in which we can cheat. Let's call a rectangle a {\em near-square} if its sides are consecutive natural numbers. Then the special Pythagorean triples in which we are interested  correspond to  the sum of two identical near-squares being another near-square.  For example, the near-square sums that correspond to $3^2+4^2=5^2$ and $20^2+21^2=29^2$ are
 \newpage

\begin{figure}[h]
\centerline{\includegraphics[scale=.5]{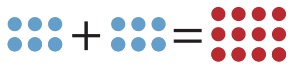}}
\centerline{$(2\times 3)+(2\times 3)=3\times 4$}
\centerline{\includegraphics[scale=.35]{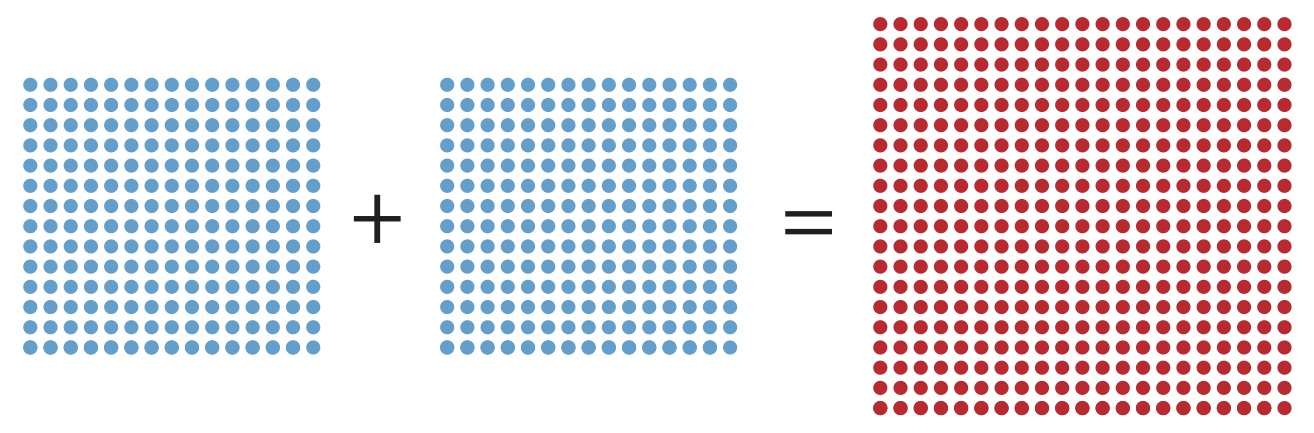}}
\centerline{$(14\times 15)+(14\times 15)=20\times 21$}\end{figure}

Now let's see what happens when we perform the merging manoeuvre on our $20 \times 21$ near-square:

\begin{figure}[h]
\centerline{\includegraphics[scale=.35]{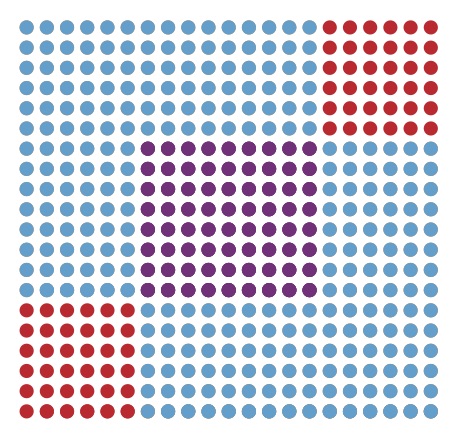}}
\centerline{Conclusion: red square + red square = purple near-square}\centerline{$6^2+6^2=8\times 9$} 
\end{figure}

So, yet another way to cheat has materialized: we have two squares adding up to a near square. Let's continue, applying the merging manoeuvre   to the two little red squares and the purple near-square:

\newpage

\begin{figure}[h]
\centerline{\includegraphics[scale=.35]{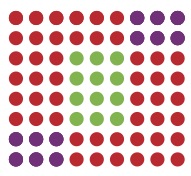}}
\centerline{Conclusion: purple near-square + purple near-square = green near-square}
\centerline{$(2\times 3)+(2\times 3)=3\times 4$}
\end{figure}

So, we've ended up with our original, tiny near-squares. 

It  is easy to check that this will always work:  beginning with two identical near-squares adding to a near-square,
and   applying the merging manoeuvre twice, we'll create two smaller  identical near-squares adding to a near-square. 
It is also not hard to see that if we continue then, no matter the near-squares with which we began, we'll wind up with the smallest 
 $(2\times 3)+(2\times 3)=3\times 4$ example. (Applying the manoeuvre again gives the $1 \times 2$ near-square, which
 gets removed entirely at the next stage.)

By reversing the procedure we can see that all such near-square sums arise from our smallest example. To illustrate, the next bigger example
to follow
$(14\times 15)+(14\times 15)=20\times 21$, and one truly worthy of the Glorious Leader, is indicated by the following two diagrams:
\begin{figure}[h]
\centerline{\includegraphics[scale=.35]{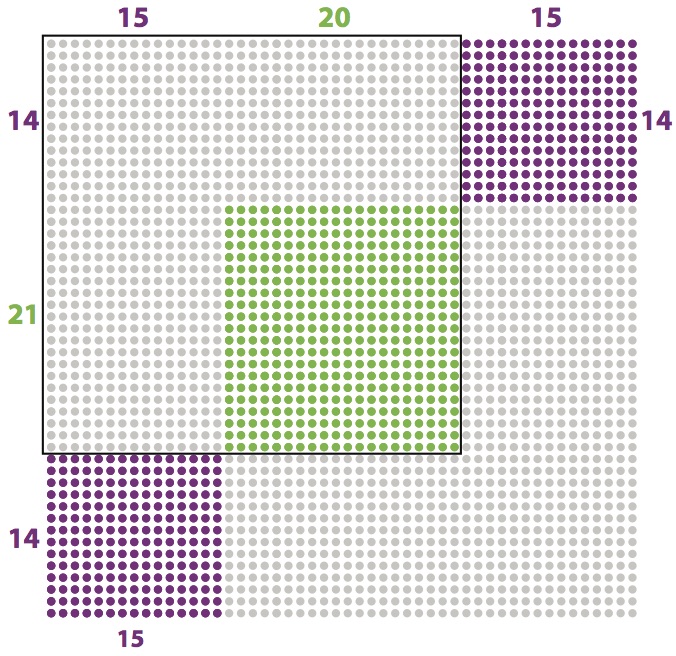}}
\end{figure}

The highlighted square has sidelength $21+14=35$, and the sidelengths of the encompassing near-square are $14+21+14=49$ and $15+20+15=50$. So, we have
$$35^2+35^2=49\times 50\,.$$ 

\newpage

Reversing the procedure again,  we obtain an impressively large near-square sum.

\begin{figure}[h]
\centerline{\includegraphics[scale=.45]{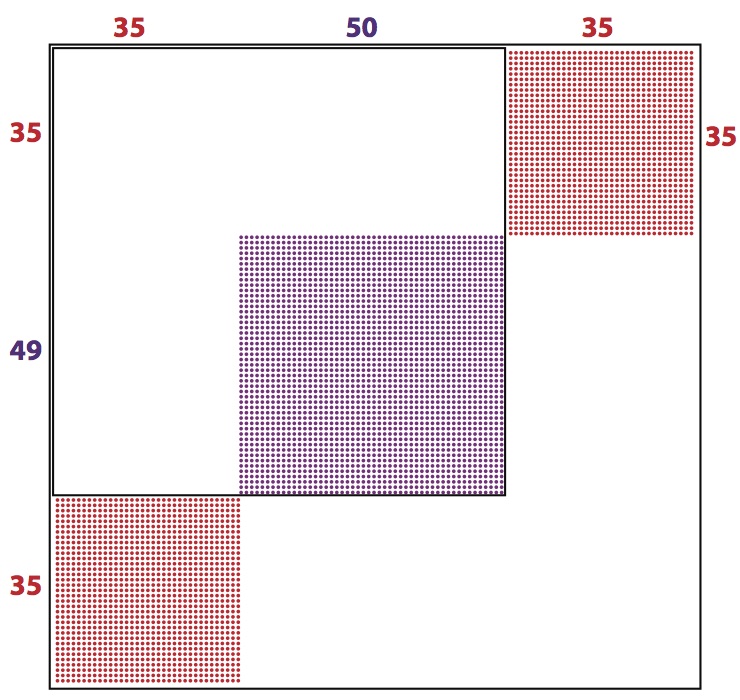}}
$$(84\times 85) + (84 \times 85)= 119 \times 120\,.$$
\end{figure}

Similar considerations allow us to  calculate the smallest examples  of our three different ways of cheating:

\begin{enumerate}
\item {\bf (Pythagorean) A square plus next-square equals a square}: $m^2+(m+1)^2=n^2$ 
$$m/n=3/5, 20/29, 119/169, 696/985,4059/5741,23660/33461,  \ldots$$
\item {\bf Twice a near-square equals a near-square}: $m(m+1)+m(m+1)=n(n+1)$
$$m/n= 2/3, 14/20, 84/119, 492/696, 2870/4059, 16730/23660, \ldots$$
\item {\bf Twice a square equals a  near-square}: $m^2+m^2=n(n+1)$
$$m/n=1/1, 6/8, 35/49, 204/288, 1189/1681, 6930/9800, 40391/57121,\ldots$$
\end{enumerate}

There are a number of other ways to cheat, including with the Pythagorean triples themselves.
In the following diagram we've begun with the equation $20^2+21^2=29^2$ and we find

\begin{figure}[h]
\centerline{\includegraphics[scale=.5]{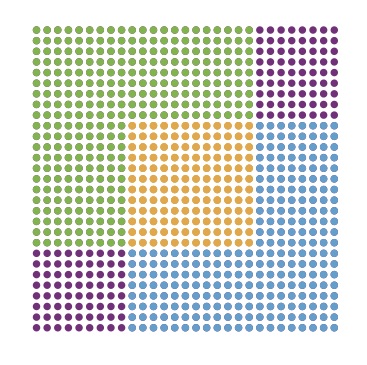}}
\centerline{purple near-square + purple near-square = yellow square}
\centerline{$(8\times 9)+(8\times 9)=12\times 12$}
\end{figure}

\newpage 
So, a fourth type of cheat is

\begin{enumerate}
\item[4.] {\bf Twice a near-square equals a square: $m(m+1)+m(m+1)=n^2$}
$$m/n=1/2,8/12,49/70,288/408,1681/2378, 9800/13860  \ldots$$
\end{enumerate}

Applying the rearranging manoeuvre to this fourth cheat produces a new Pythagorean triple. Starting with the $(8 \times 9) + (8 \times 9) = 12 \times 12$ example, we obtain our familiar friend $3^2 + 4^2 = 5^2$.

\begin{figure}[h]
\centerline{\includegraphics[scale=.5]{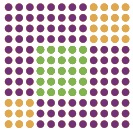}}
\centerline{ \rm yellow square +  yellow square = green square}
\centerline{$3^2+4^2=5^2$}
\end{figure}

There is one more cheat-variation worth noting. Notice that there are two different ways of arranging  two near-squares into a larger square, differing by  the near-squares' relative orientation. We've employed one orientation in our fourth cheat above, and the other orientation results in a failed cheat, as illustrated below. 
\newpage

\begin{figure}[h]
\centerline{\includegraphics[scale=.5]{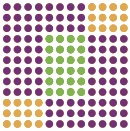}}
\centerline{yellow near-square + yellow near-square = green not-near-square}
\centerline{$(3\times 4)+(3\times 4)=4\times 6$}
\end{figure}

That wasn't very successful but something new occurs when we reconsider our original scenario, fitting two identical near-squares into a larger near-square. There are in effect three different possible orientations.
 
\begin{figure}[h]
\centerline{\includegraphics[scale=.35]{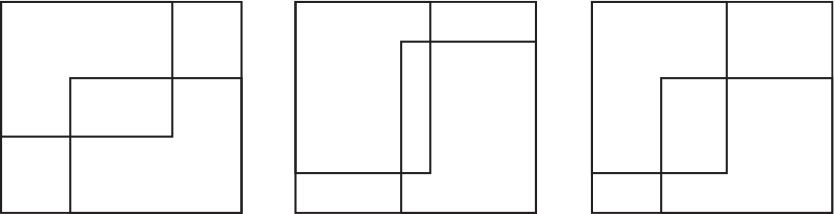}}
\end{figure}
We used the first orientation in our second cheat above, and the second orientation produces not-square-or-near-squares.  However, the third orientation gives us a fifth and final method  of cheating. Below we've merged two $14\times 15$ near-squares inside a $20\times 21$ in this way:

\begin{figure}[h]
\centerline{\includegraphics[scale=.35]{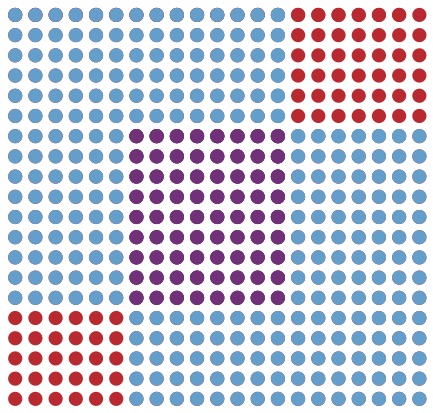}}
\centerline{Conclusion:  red near-square + red next-near-square =   purple near-square}
\centerline{$(5\times 6)+(6\times 7)=8\times 9$}
\end{figure}

In general, this cheat gives us

\begin{enumerate}
\item[5.] {\bf  near-square plus next-near-square equals near-square: 

$m(m+1)+(m+1)(m+2)=n(n+1) $ }
$$m/n=5/8,34/49,203/288,1188/1681, 6929/9800, \ldots  $$
\end{enumerate}

The general solutions to our five cheat equations are given by recursion formulas, and it is straight-forward to derive these algebraic formulas from our geometric merging procedures.  However the truth
is that   everything about integer solutions to these five equations has been known for a long time. In particular, the special Pythagorean triples we've been pondering form one type of the so-called {\em Pythagorean twins}. General formulas can be found in many places, including at the {\em Online Encyclopedia of Integer Sequences} \cite{EIS}.

\section{Impossible cheating}

At this point it also makes sense to list all potential cheating schemes of a reasonably similar nature, and to figure out which work and which do not. We  should then consider 
which of  the following eight equations have natural numbers solutions:
\begin{enumerate}[(a)]
\item $m^2+m^2=n^2$: Impossible, by Tennenbaum's argument.
\item $m^2+m^2=n(n+1)$: Possible, our 3rd cheat.
\item $m^2+(m+1)^2=n^2$: Possible, Pythagorean twins, our 1st cheat.
\item $m^2+(m+1)^2=n(n+1)$: Not possible; see below.
\item $m(m+1)+m(m+1)=n^2$: Possible, our 4th cheat.
\item $m(m+1)+m(m+1)=n(n+1)$: Possible, our 2nd cheat.
\item $m(m+1)+(m+1)(m+2)=n^2$: Not possible; see below.
\item $m(m+1)+(m+1)(m+2)=n(n+1)$: Possible, our 5th cheat.
\end{enumerate}

We'll now show that the final cases, (d) and (g), are impossible. 

First to (d). It is very easy to see that  the equation $m^2+(m+1)^2=n(n+1)$  is impossible since the left side  would be odd and the right side would be even.

For (g), the equation $m(m+1)+(m+1)(m+2)=n^2$ can be rewritten as $2 (m+1)^2 = n^2$, which reduces to case (a) above, and the irrationality of $\sqrt 2$. 

What is interesting is that the impossibility of (g) can be proved with the same merging proof by contradiction with which we began the article. This is illustrated below.
\newpage

\begin{figure}[h]
\centerline{\includegraphics[scale=.45]{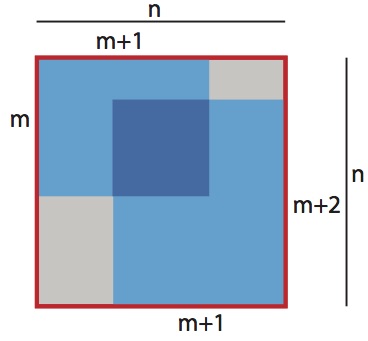}}
\end{figure}
It is easy to check that the two grey rectangles are a near-square and next-near-square, and that the  overlapped region is a square. It follows that
 any supposed minimal example of $m(m+1)+(m+1)(m+2)=n^2$ being satisfied would result an even smaller example, which is impossible. 

Since, as we noted,  $m(m+1)+(m+1)(m+2)=n^2$ can be rewritten as $2 (m+1)^2 = n^2$, this then  amounts to an alternative proof that $\sqrt 2$ is irrational.

\section{Bending the rules without cheating} 

Recall that the {\em triangular numbers} are the numbers of the form $\frac{n(n+1)}{2}$. They take their name from counting the number of objects  in an  equilateral triangle. 
Of course, this also means that any of our near-square identities, such as 
$$(14\times 15)+(14\times 15)=20\times 21\,,$$
represents two  equilateral triangles merging into a larger  equilateral triangle.
\newpage

\begin{figure}[h]
\centerline{\includegraphics[scale=.35]{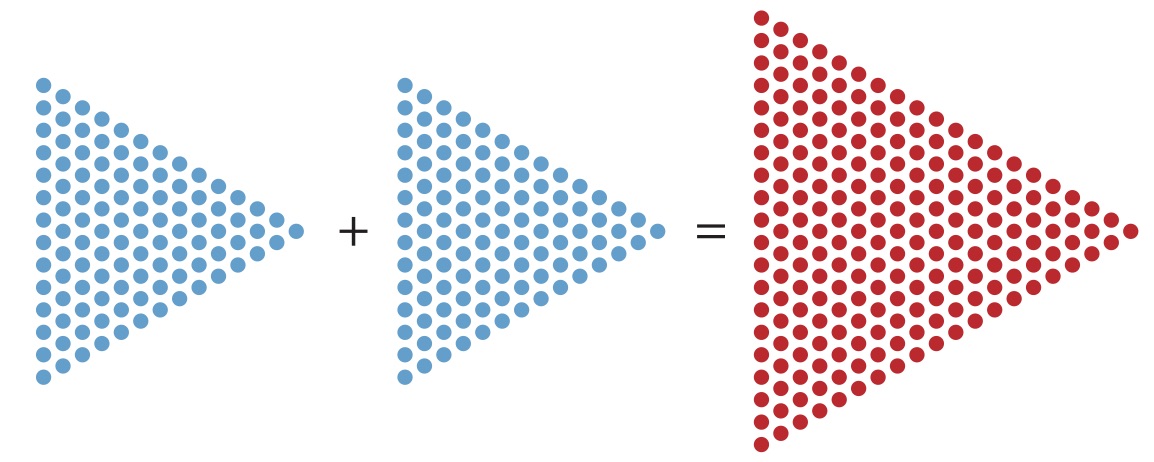}}
$$\frac{14\times 15}{2}+\frac{14\times 15}{2}=\frac{20\times 21}{2}$$
\end{figure}

Maybe our glorious leader also likes equilateral triangles (in which case we'll show him this: \url{http://youtu.be/rjHFkx6eTL8})? But if not, and if it's got to be squares, perhaps a good way to sell the Pythagorean twin cheat is to show that a square and next-square  can be merged into a square in two different and very interesting ways:

\begin{figure}[h]
\centerline{\includegraphics[scale=.35]{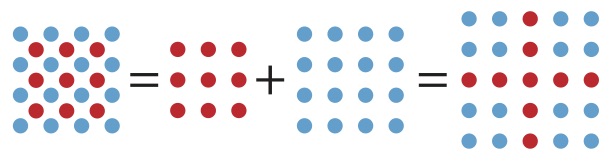}}
$$3^2+4^2=5^2$$
\end{figure}

\section{A little perspective on marching squares} This article contains much that is very well known. So what's the point? 

The   first highlight of our story is the clever covering argument in Tennenbaum's proof that $\sqrt 2$ is an irrational number. We believe  that this very beautiful proof deserves to be more widely known. However, there is a problem.

 Tennenbaum's is a proof by contradiction and people with little mathematical background, especially young students being led into the wonders of mathematics,  almost invariably  struggle with the basic line of reasoning in such  proofs.
We have attempted to present  Tennenbaum's lovely proof so that the contradiction aspect is not even noticable until the final summing up, when we announce  that a proof by contradiction is what we've just done. We've had pleasing success with this story approach in public lectures, school and university presentations, and even in a recent column in a newspaper here in Australia; see  \cite{age}. Of course selling tricky mathematics with stories is hardly new. There is even an entire, very good book about $\sqrt 2$ that treads similar ground; see \cite{flannery}. 

In the second part of the article we do something different, and apparently new. We take the ``destructive'' argument within a proof by contradiction and turn it around to construct examples of interesting and very pretty near-misses. It may be fruitful to  explore other  proofs by contradiction, to see whether constructions can be flipped in this surprising manner. Conway and Shipman \cite{conway1} mention  other covering arguments which are similar in spirit to Tennenbaum's proof; these and  many other ``infinite descent'' geometric proofs may bear fruit. 

In the third part of the article we present an interesting variation of Tennenbaum's proof of the irrationality of $\sqrt 2$.

Finally, if somewhat optimistically, we hope that we've provided some food for thought for anyone in the business of choreographing marching bands and parades.

\end{document}